# GLOBAL STABILITY RESULTS FOR SYSTEMS UNDER SAMPLED-DATA CONTROL

**Iasson Karafyllis**[*] **and Costas Kravaris**[**]

[*]**Department of Environmental Engineering, Technical University of Crete, 73100, Chania, Greece**
email: ikarafyl@enveng.tuc.gr

[**]**Department of Chemical Engineering, University of Patras
1 Karatheodory Str., 26500 Patras, Greece**
email: kravaris@chemeng.upatras.gr

**Abstract**
In this work sufficient conditions expressed by means of single and vector Lyapunov functions of Uniform Input-to-Output Stability (UIOS) and Uniform Input-to-State Stability (UISS) are given for finite-dimensional systems under feedback control with zero order hold.

## 1. Introduction

Given the finite-dimensional continuous-time system:

$$\begin{aligned} \dot{x}(t) &= f(x(t), u(t), d(t)) \\ Y(t) &= H(x(t)) \\ x(t) &\in \Re^n, u(t) \in U \subseteq \Re^m, d(t) \in D \subseteq \Re^l, Y(t) \in \Re^k \end{aligned} \quad (1.1)$$

where the vector field $f: \Re^n \times U \times D \to \Re^n$ is continuous, locally Lipschitz in $x \in \Re^n$, $u(t)$ represents the control input and $d(t)$ unknown disturbances or model uncertainty. Suppose now that a continuous-time control law $u(t) = k(x(t))$ is applied to the system (1.1). Then the resulting closed-loop system can be analyzed and the feedback function can be selected for desirable properties, e.g., robust global asymptotic stability or Input-to-State Stability.

When the foregoing control law is implemented in discrete time under zero-order hold with sampling period $r$

$$\begin{aligned} u(t) &= k(x(\tau_i)), t \in [\tau_i, \tau_{i+1}) \\ \tau_{i+1} &= \tau_i + r, i = 0,1,2,... \end{aligned} \quad (1.2)$$

the resulting closed-loop system

$$\begin{aligned} \dot{x}(t) &= f(x(t), k(x(\tau_i)), d(t)), t \in [\tau_i, \tau_{i+1}) \\ \tau_{i+1} &= \tau_i + r, i = 0,1,2,... \\ Y(t) &= H(x(t)) \end{aligned} \quad (1.3)$$

does not necessarily possess the same properties as the closed-loop system under continuous-time implementation. Also, the deduction of properties under discrete-time feedback from properties under continuous-time feedback is quite difficult and involved. Questions like

"Can the state feedback $u = k(x)$ robustly globally asymptotically stabilize the origin under zero-order hold discrete-time implementation, if it globally asymptotically stabilizes the origin under continuous-time implementation?"

do not have an affirmative answer in general, even if the sampling period $r$ is assumed to be arbitrarily small. However, if we can guarantee an affirmative answer to the previous question then, the next important question that arises is:

"What is the maximum allowable sampling period for guaranteed global asymptotic stability?"

In general, the task of providing answers to the above questions is highly non-trivial.

On the other hand, the analysis of the closed-loop system (1.3) as a discrete-time system (ignoring inter-sample behavior) is possible in principle, but it is hindered by the unavailability (in general) of an exact sampled-data representation of the given continuous-time nonlinear system.

There is a large body of literature concerning the foregoing very important and very challenging issues. In particular, the following lines of attack have been pursued to derive stability results:

* making use of numerical approximations of the solution of the open-loop system (e.g., in the work of D. Nesic, A. Teel and others, see [10,11,20,24,25,31-37,49]). The results obtained in this way lead to a systematic procedure for the construction of practical, semi-global feedback stabilizers and provide a list of possible reasons that explain the occasional failure of sampled-data control mechanisms. Recent research takes into account performance and robustness issues as well (see [10,21,23,37]).
* exploiting special characteristics of the system such as homogeneity (see [8]), global Lipschitz conditions (see [13]) or linear structure with uncertainties (see [2] as well as the textbook [42]).
* making use of Linear Matrix Inequalities in the context of hybrid systems (see [14,15,27,29,48]).
* considering the closed-loop system as a discrete-time system, i.e. ignoring inter-sample behavior (see [1,28]). As mentioned earlier, this point of view generally requires the knowledge of the solution map, which is rarely available for the nonlinear case. However, recent work has established results that characterize the inter-sample behavior of the solutions based on the behavior of the solution of the discrete-time system (see [31]).

It is important to note that the above very important research results do not provide conditions for global Asymptotic Stability or Input-to-Output Stability for general nonlinear sampled-data systems (usually only semi-global practical stability properties are established or global stability for limited classes of systems).

A very interesting point of view that was recently explored in the linear systems theory (see [5,6]) involves considering the closed-loop system under zero order-hold as a time-delay system. This is a natural and intuitively meaningful point of view, since sampled-data control introduces piecewise constant retarded arguments. It is exactly this point of view that will be explored in the present work for nonlinear systems for the first time.

Moreover, the present work will utilize the method of Razumikhin functions (see [12,38,47]) for stability analysis of time-delay systems, together with recent developments in the theory of vector Lyapunov functions (see [30]). The main contribution will be the development of Lyapunov-like sufficient conditions for Uniform Input-to-Output Stability (UIOS) and Uniform Input-to-State Stability (UISS). The notions of UIOS and UISS were formulated in [39,40,41,44,45] for finite-dimensional systems described by ordinary differential equations. More recently, sufficient conditions for semiglobal practical UISS were studied in [24,25] for sampled-data systems. In the present work, the sufficient conditions will be expressed in terms of a scalar and in terms of a vector Lyapunov function. The main features of the results will be that they

* can be applied to uncertain nonlinear systems with no special characteristics (such as homogeneity or global Lipschitz conditions),
* can lead to global stability results for an arbitrary output map $Y = H(x)$ (including the possibility of using the identity map $H(x) \equiv x$ and thus recovering the familiar stability notions for the state)
* do not require knowledge of the solution map for the open-loop system,
* allow the direct determination of the maximum allowable sampling period,
* guarantee robustness to perturbations of the sampling schedule

The present work will not aim directly at the construction/design of sampled-data feedback stabilizers. The main goal will be to develop, for the first time, Lyapunov-like sufficient conditions for robust global asymptotic stability or Input-to-Output Stability of (general nonlinear) sampled-data systems. These Lyapunov-like sufficient conditions in the present work can be potentially useful for the construction and/or the Lyapunov redesign of sampled-data stabilizers for nonlinear control systems, although this application will not be addressed in detail in the present work.

More specifically, in this work we will consider hybrid systems which are described in the following way: given a pair of sets $D \subseteq \Re^l$, a positive function $h : \Re^n \to (0, r]$, which is bounded by a certain constant $r > 0$ and a pair of vector fields $f : \Re^n \times \Re^n \times D \times U \times U \to \Re^n$, $H : \Re^n \to \Re^p$, we consider the hybrid system that produces for each $(t_0, x_0) \in \Re^+ \times \Re^n$ and for each triplet of measurable and locally bounded inputs $d : \Re^+ \to D$, $\tilde{d} : \Re^+ \to \Re^+$, $v : \Re^+ \to U$ the piecewise absolutely continuous function $t \to x(t) \in \Re^n$, produced by the following algorithm:

Step $i$:
1) Given $\tau_i$ and $x(\tau_i)$, calculate $\tau_{i+1}$ using the equation $\tau_{i+1} = \tau_i + \exp(-\tilde{d}(\tau_i))h(x(\tau_i))$,
2) Compute the state trajectory $x(t)$, $t \in [\tau_i, \tau_{i+1})$ as the solution of the differential equation $\dot{x}(t) = f(x(t), x(\tau_i), d(t), v(t), v(\tau_i))$,
3) Calculate $x(\tau_{i+1})$ using the equation $x(\tau_{i+1}) = \lim_{t \to \tau_{i+1}^-} x(t)$.

For $i = 0$ we take $\tau_0 = t_0$ and $x(\tau_0) = x_0$ (initial condition). Schematically, we write

$$\begin{aligned} &\dot{x}(t) = f(x(t), x(\tau_i), d(t), v(t), v(\tau_i)) \quad, \quad t \in [\tau_i, \tau_{i+1}) \\ &\tau_0 = t_0, \tau_{i+1} = \tau_i + \exp(-\tilde{d}(\tau_i))h(x(\tau_i)), i = 0, 1, ... \\ &Y(t) = H(x(t)) \end{aligned} \qquad (1.4)$$

with initial condition $x(t_0) = x_0$. Hybrid systems of the form (1.4) will be called "sampled-data" systems.

Systems of the form (1.4) arise as closed-loop finite-dimensional systems under feedback control with zero order hold. Specifically, considering the finite-dimensional continuous-time dynamic system

$$\begin{aligned} &\dot{x}(t) = f(x(t), u(t), d(t)) \\ &Y(t) = H(x(t)) \end{aligned}, \quad x(t) \in \Re^n, u(t) \in \Re^m, d(t) \in D \qquad (1.5)$$

where the vector field $f : \Re^n \times \Re^m \times D \to \Re^n$ is continuous, locally Lipschitz in $x \in \Re^n$. Then, the application of the (not necessarily continuous) feedback law

$$u(t) = k(x(\tau_i) + e(\tau_i)) + v(t) \text{ on the interval } [\tau_i, \tau_i + \exp(-\tilde{d}(\tau_i))h(x(\tau_i))), \ i = 0, 1, 2, ... \qquad (1.6)$$

where $v(t) \in U \subseteq \Re^m$ represents the control actuator error, $e(t) \in \Re^n$ represents the measurement error and $h : \Re^n \to (0, r]$ is a positive function bounded by certain constant $r > 0$ and $\tilde{d}(t) \geq 0$ is the input that quantifies the perturbation of the sampling schedule, gives rise to the dynamical system, that produces for each $(t_0, x_0) \in \Re^+ \times \Re^n$ and for each measurable and locally bounded inputs $d : \Re^+ \to D$, $v : \Re^+ \to U$, $e : \Re^+ \to \Re^n$ and $\tilde{d} : \Re^+ \to \Re^+$ the absolutely continuous function $[t_0, +\infty) \ni t \to x(t) \in \Re^n$ that satisfies a.e. the differential equation

$$\begin{aligned} &\dot{x}(t) = f(x(t), k(x(\tau_i) + e(\tau_i)) + v(t), d(t)) \quad, \quad t \in [\tau_i, \tau_{i+1}) \\ &\tau_0 = t_0, \tau_{i+1} = \tau_i + \exp(-\tilde{d}(\tau_i))h(x(\tau_i)), i = 0, 1, ... \\ &Y(t) = H(x(t)) \end{aligned} \qquad (1.7)$$

with initial condition $x(t_0) = x_0$. Clearly, system (1.7) is a system of the form (1.4).

Some comments are in order relative to the discrete-time feedback law (1.6). In addition to accounting for actuator and measurement errors, (1.6) is a generalization of (1.2) in the following ways:

i) The sampling period is allowed to be state-dependent. The situation of a non-constant $h: \Re^n \to (0, r]$ has recently emerged in connection with the study of asymptotic controllability of nonlinear systems (see the classical work in [3,43] and the links with sampled-data stabilizability given in [7]). Of course, the situation of $h(x) \equiv r > 0$ (constant) is the one of practical significance, in which case we say that the sampled-data feedback law $u(t) = k(x(\tau_i) + e(\tau_i)) + v(t)$ is applied with positive sampling rate (see [7]).

ii) The sampling period is allowed to be time-varying. The factor $\exp(-\tilde{d}(\tau_i)) \leq 1$, with $\tilde{d}(t) \geq 0$ some non-negative function of time, is an uncertainty factor in the end-point of the sampling interval. Proving stability for any non-negative input $\tilde{d}: \Re^+ \to \Re^+$ will guarantee stability for all sampling schedules with $\tau_{i+1} - \tau_i \leq h(x(\tau_i))$ (robustness to perturbations of the sampling schedule).

The structure of the present paper is as follows: In Section 2, we state our main assumptions for system (1.4). Using the system-theoretic framework presented in [18,19], we exploit the results of [17] and demonstrate that system (1.4) under the hypotheses of Section 2 satisfies essential properties for its qualitative study. In Section 3, we state our main results. In Section 4, illustrating examples are presented, which show the flexibility that a vector Lyapunov function can provide. Moreover, we show how the results of Section 3 can be used in conjunction with the backstepping method for triangular nonlinear systems to check whether the feedback, constructed by the backstepping method presented in [4], can robustly globally asymptotically stabilize the equilibrium point when applied with zero order hold and positive sampling rate. Finally, Section 5 contains the concluding remarks of the present work.

**Notations** Throughout this paper we adopt the following notations:

* Let $I \subseteq \Re$ be an interval. By $C^0(I; \Omega)$, we denote the class of continuous functions on $I$, which take values in $\Omega$. By $C^1(I; \Omega)$, we denote the class of functions on $I$ with continuous derivative, which take values in $\Omega$.

* For a vector $x \in \Re^n$ we denote by $|x|$ its usual Euclidean norm and by $x'$ its transpose.

* We denote by $K^+$ the class of positive $C^0$ functions defined on $\Re^+$. We say that a non-decreasing continuous function $\gamma : \Re^+ \to \Re^+$ is of class $N$ if $\gamma(0) = 0$. We say that a function $\rho : \Re^+ \to \Re^+$ is positive definite if $\rho(0) = 0$ and $\rho(s) > 0$ for all $s > 0$. For the definitions of the classes $K$ and $K_\infty$, see [22]. By $KL$ we denote the set of all continuous functions $\sigma = \sigma(s,t) : \Re^+ \times \Re^+ \to \Re^+$ with the properties: (i) for each $t \geq 0$ the mapping $\sigma(\cdot, t)$ is of class $K$; (ii) for each $s \geq 0$, the mapping $\sigma(s, \cdot)$ is non-increasing with $\lim_{t \to +\infty} \sigma(s,t) = 0$.

* Let $U \subseteq \Re^m$ be a non-empty set with $0 \in U$. By $B_U[0,r] := \{ u \in U ; |u| \leq r \}$ we denote the closed sphere in $U \subseteq \Re^m$ with radius $r \geq 0$, centered at $0 \in U$.

* Let $D \subseteq \Re^l$ be a non-empty set. By $M(D)$ we denote the class of all Lebesgue measurable and locally bounded mappings $d : \Re^+ \to D$.

* For every scalar function $V \in C^1(\Re^n; \Re)$, $\nabla V(x)$ denotes the gradient of $V \in C^1(\Re^n; \Re)$ at $x \in \Re^n$, i.e., $\nabla V(x) = \left( \frac{\partial V}{\partial x_1}(x), \ldots, \frac{\partial V}{\partial x_n}(x) \right)$. We say that a function $V \in C^0(\Re^n; \Re)$ is positive definite if $V(x) > 0$ for all $x \neq 0$ and $V(0) = 0$. We say that a function $V \in C^0(\Re^n; \Re)$ is radially unbounded if the following property holds: "$V(x)$ is bounded if and only if $|x|$ is bounded".

## 2. Main Assumptions and Preliminaries for Sampled-Data Systems

In the present work we study systems of the form (1.4) under the following hypotheses:

**(H1)** $f(x, x_0, d, v, v_0)$ is continuous with respect to $(x, d, v) \in D \times U$ and such that for every bounded $S \subset \Re^n \times \Re^n \times U \times U$ there exists constant $L \geq 0$ such that

$$(x-y)'(f(x, x_0, d, v, v_0) - f(y, x_0, d, v, v_0)) \leq L|x-y|^2$$
$$\forall (x, x_0, v, v_0, d) \in S \times D, \forall (y, x_0, v, v_0, d) \in S \times D \quad (2.1)$$

**(H2)** There exist a function $a \in K_\infty$ such that

$$|f(x, x_0, d, v, v_0)| \leq a(|x| + |x_0| + |v| + |v_0|), \quad \forall (v, v_0, d, x, x_0) \in U \times U \times D \times \Re^n \times \Re^n \quad (2.2)$$

**(H3)** $H : \Re^n \to \Re^p$ is a continuous map with $H(0) = 0$. Moreover, there exists a constant $R \geq 0$ and a function $p \in K_\infty$ such that $|x| \leq R + p(|H(x)|)$.

**(H4)** The function $h : \Re^n \to (0, r]$ is a positive, continuous and bounded function.

The following theorem is an immediate consequence of the results presented in [18].

**Theorem 2.1:** *Consider the control system (1.4) under hypotheses (H1-4) and let $\phi(t, t_0, x_0; v, d, \tilde{d})$ denote the solution of (1.4) with initial condition $x(t_0) = x_0$ corresponding to input $(d, \tilde{d}, v) \in M(D) \times M(\Re^+) \times M(U)$. Then*

(i) *system (1.4) has the **"Boundedness-Implies-Continuation"** (BIC) property, i.e., for each $(t_0, x_0, v, d, \tilde{d}) \in \Re^+ \times \Re^n \times M(U) \times M(D) \times M(\Re^+)$, there exists $t_{max} \in (t_0, +\infty]$ (the maximal existence time of the solution) such that the solution $\phi(t, t_0, x_0; v, d, \tilde{d})$ of (1.4) exists for all $t \in [t_0, t_{max})$. In addition, if $t_{max} < +\infty$ then for every $M > 0$ there exists $t \in [t_0, t_{max})$ with $|\phi(t, t_0, x_0; v, d, \tilde{d})| > M$.*

(ii) *$0 \in \Re^n$ is a **robust equilibrium point** from the input $(v, \tilde{d}) \in M(U) \times M(\Re^+)$, i.e., for every $\varepsilon > 0$, $T \in \Re^+$ there exists $\delta := \delta(\varepsilon, T) > 0$ such that for all $(t_0, x_0, v, \tilde{d}) \in [0, T] \times \Re^n \times M(U) \times M(\Re^+)$, with $|x_0| + \sup_{t \geq 0} |v(t)| + \sup_{t \geq 0} \tilde{d}(t) < \delta$ it holds that the solution $\phi(t, t_0, x_0; v, d, \tilde{d})$ of (1.4) exists for all $t \in [t_0, t_0 + T]$ and $d \in M(D)$ and*

$$\sup\{|\phi(t, t_0, x_0; v, d, \tilde{d})|; d \in M(D), t \in [t_0, t_0 + T], t_0 \in [0, T]\} < \varepsilon$$

(iii) *system (1.4) is **autonomous**, i.e., for each $(t_0, x_0, v, d, \tilde{d}) \in \Re^+ \times \Re^n \times M(U) \times M(D) \times M(\Re^+)$, $t \geq t_0$ and for each $\theta \in (-\infty, t_0]$ it holds that $\phi(t, t_0, x_0, v, d, \tilde{d}) = \phi(t - \theta, t_0 - \theta, x_0, P_\theta v, P_\theta d, P_\theta \tilde{d})$, where $(P_\theta v)(t) = v(t + \theta)$, $(P_\theta \tilde{d})(t) = \tilde{d}(t + \theta)$ and $(P_\theta d)(t) = d(t + \theta)$ for all $t + \theta \geq 0$.*

We next present the notion of Uniform Input-to-Output Stability property for systems of the form (1.4) under hypotheses (H1-4). The following definition is parallel to the corresponding notions used for finite-dimensional control systems described by ordinary differential equations (see [39,40,41,44,45]).

**Definition 2.2:** *Consider the control system (1.4) under hypotheses (H1-4) and let $\phi(t, x_0; v, d, \tilde{d})$ denote the solution of (1.4) with initial condition $x(0) = x_0$ corresponding to input $(d, \tilde{d}, v) \in M(D) \times M(\Re^+) \times M(U)$. Suppose that (1.4) is **robustly forward complete** (RFC) from the input $(v, \tilde{d}) \in M(U) \times M(\Re^+)$, i.e., suppose that for every $R \geq 0$, $T \geq 0$, it holds that*

$$\sup\{|\phi(t, x_0; v, d, \tilde{d})|; v \in M(B_U[0, R]), \tilde{d} \in M([0, R]), t \in [0, T], |x_0| \leq R, d \in M(D)\} < +\infty \quad (2.3)$$

*Moreover, suppose that there exist functions $\sigma \in KL$, $\gamma \in N$ such that the following estimate holds for all $(x_0, v, d, \tilde{d}) \in \Re^n \times M(U) \times M(D) \times M(\Re^+)$ and $t \geq 0$:*

$$|H(\phi(t, x_0; v, d, \tilde{d}))| \leq \sigma(|x_0|, t) + \sup_{0 \leq \tau \leq t} \gamma(|v(\tau)|) \quad (2.4)$$

*Then we say that (1.4) satisfies the **Uniform Input-to-Output Stability (UIOS) property** with gain $\gamma \in N$ from the input $v \in M(U)$ and zero gain from the input $\tilde{d} \in M(\Re^+)$. Particularly, if $H(x) := x$, then we say that (1.4) satisfies the Uniform Input-to-State Stability (UISS) property with gain $\gamma \in N$ from the input $v \in M(U)$ and zero gain from the input $\tilde{d} \in M(\Re^+)$.*

For the proof of our main results we will need the following technical small-gain lemma. It is a direct corollary of Theorem 1 in [46] and is closely related to Lemma A.1 in [16].

**Lemma 2.3:** *For every $\sigma \in KL$ and $a \in K$ with $a(s) < s$ for all $s > 0$, there exists $\tilde{\sigma} \in KL$ with the following property: if $y : [t_0, t_1) \to \Re^+$, $u : \Re^+ \to \Re^+$ are locally bounded functions and $M \geq 0$ a constant such that the following inequality holds for all $\xi \in [t_0, t_1)$:*

$$y(t) \leq \max\left\{ \sigma(M, t-\xi) \; ; \; a\left(\sup_{\xi \leq \tau \leq t} y(\tau)\right) ; \; u(t) \right\}, \; \forall t \in [\xi, t_1) \tag{2.5}$$

*then the following estimate holds for all $t \in [t_0, t_1)$:*

$$y(t) \leq \max\left\{ \tilde{\sigma}(M, t-t_0) \; ; \; \sup_{t_0 \leq \tau \leq t} u(\tau) \right\}, \; \forall t \in [t_0, t_1) \tag{2.6}$$

Finally, we end this section by presenting the following comparison lemma, which provides a sharp estimate of the evolution of Lyapunov functions (compare the obtained estimate with Theorem 5.2, page 218 in [22]). Its proof can be found in the Appendix of the present work. Notice that a similar result with the following lemma is contained in [9].

**Lemma 2.4:** *For each positive definite continuous function $\rho : \Re^+ \to \Re^+$ there exists a function $\sigma$ of class $KL$, with $\sigma(s,0) = s$ for all $s \geq 0$ with the following property: if $y : [t_0, t_1] \to \Re^+$ is an absolutely continuous function, $u : \Re^+ \to \Re^+$ is a locally bounded mapping and $I \subset [t_0, t_1]$ a set of Lebesgue measure zero such that $\dot{y}(t)$ is defined on $[t_0, t_1] \setminus I$ and such that the following implication holds for all $t \in [t_0, t_1] \setminus I$:*

$$y(t) \geq u(t) \;\; \Rightarrow \;\; \dot{y}(t) \leq -\rho(y(t)) \tag{2.7}$$

*then the following estimate holds for all $t \in [t_0, t_1]$:*

$$y(t) \leq \max\left\{ \sigma(y(t_0), t-t_0), \sup_{t_0 \leq s \leq t} \sigma(u(s), t-s) \right\} \tag{2.8}$$

## 3. Main Results and Examples

The main result of the present work is the following theorem, which provides sufficient conditions for the UIOS property for system (1.4) with gain $\gamma \in N$ from the input $u \in M(U)$ and zero gain from the input $\tilde{d} \in M(\Re^+)$. The conditions are expressed by means of a vector Lyapunov function. It utilizes the method of Razumikhin functions (see [12,38,47]) for stability analysis of time-delay systems, together with recent developments in the theory of vector Lyapunov functions (see [30]).

**Theorem 3.1 (Vector Lyapunov Function Characterization of UIOS):** *Consider system (1.4) under hypotheses (H1-4) and suppose that there exists a family of functions $V_i \in C^1(\Re^n; \Re^+)$ ($i = 1,...,k$), functions $a_1, a_2 \in K_\infty$ $a, \zeta \in N$ with $a(s) < s$ for all $s > 0$ and a family of positive definite functions $\rho_i \in C^0(\Re^+; \Re^+)$ ($i = 1,...,k$), such that the following inequalities hold:*

$$a_1(|H(x)|) \leq \max_{i=1,\ldots,k} V_i(x) \leq a_2(|x|), \quad \forall x \in \Re^n \tag{3.1}$$

$$\nabla V_i(x) f(x, x_0, d, v, v_0) \leq -\rho_i(V_i(x)),$$
$$\forall (x, d, v, v_0) \in \Re^n \times D \times U \times U \text{ with } \max\{\zeta(|v|), \zeta(|v_0|)\} \leq V_i(x) \text{ and } x_0 \in A_i(h(x_0), x) \ (i=1,\ldots,k) \tag{3.2}$$

where the family of set-valued map $\Re^+ \times \Re^n \ni (h, x) \to A_i(h, x) \subseteq \Re^n$ ($i = 1,\ldots,k$) is defined by

$$A_i(h, x) = \bigcup_{0 \leq s \leq h} \left\{ \begin{array}{l} x_0 \in \Re^n : \exists (d, v) \in M(D) \times M(U) \text{ with } \phi(s, x_0; d, v) = x, \\ \zeta(|v(t)|) \leq V_i(x), a(V_j(\phi(t, x_0; d, v))) \leq V_i(x) \\ \text{for all } t \in [0, s] \text{ and } j = 1,\ldots,k \end{array} \right\} \tag{3.3}$$

and $\phi(t, x_0; d, v)$ denotes the solution of $\dot{x}(t) = f(x(t), x_0, d(t), v(t), v(0))$ with initial condition $x(0) = x_0$ corresponding to $(d, v) \in M(D) \times M(U)$.

Then (1.4) satisfies the UIOS property with gain $\gamma = a_1^{-1} \circ \zeta \in N$ from the input $v \in M(U)$ and zero gain from the input $\tilde{d} \in M(\Re^+)$.

**Remark 3.2:** In general it is very difficult to obtain an accurate description of the set-valued map $\Re^+ \times \Re^n \ni (h, x) \to A(h, x) \subseteq \Re^n$ defined by (3.3). However, for every $g \in C^1(\Re^n; \Re)$, we have:

$$A_i(h, x) \subseteq B_i^g(h, x) = \left\{ x_0 \in \Re^n : |g(x_0) - g(x)| \leq h b_i^g(x) \right\}, \quad \forall (h, x) \in \Re^+ \times \Re^n$$

where

$$b_i^g(x) := \max\left\{ |\nabla g(\xi) f(\xi, x_0, d, v, v_0)| : d \in D, \ \zeta(\max\{|v|, |v_0|\}) \leq V_i(x), a(\max\{V(\xi), V(x_0)\}) \leq V_i(x) \right\} < +\infty$$

and $V(x) = \max_{i=1,\ldots,k} V_i(x)$. In order to justify the inclusion $A_i(h, x) \subseteq B_i^g(h, x)$, notice that for all $(d, v, x_0) \in M(D) \times M(U) \times \Re^n$ it holds that

$$g(\phi(s, x_0; d, v)) - g(x_0) = \int_0^s \nabla g(\phi(t, x_0; d, v)) f(\phi(t, x_0; d, v), x_0, d(t), v(t), v(0)) dt$$

for all $s \geq 0$ for which $\phi(s, x_0; d, v)$ exists (i.e., the solution of $\dot{x}(t) = f(x(t), x_0, d(t), v(t), v(0))$ with initial condition $x(0) = x_0$ corresponding to $(d, v) \in M(D) \times M(U)$ exists for all $t \in [0, s]$). If $x_0 \in A_i(h, x)$, by virtue of definition (3.3), there exists $s \in [0, h]$, $(d, v) \in M(D) \times M(U)$ with $\zeta(|v(t)|) \leq V_i(x)$ for all $t \in [0, s]$ such that $a(V(\phi(t, x_0; d, v))) \leq V_i(x)$ for all $t \in [0, s]$ and $\phi(s, x_0; d, v) = x$ (where $V(x) = \max_{i=1,\ldots,k} V_i(x)$). Clearly, the previous equality, implies that if $x_0 \in A_i(h, x)$ then there exists $s \in [0, h]$ such that:

$$|g(x) - g(x_0)| \leq s b_i^g(x) \leq h b_i^g(x)$$

which shows that $x_0 \in B_i^g(h, x)$.

Thus in order to establish the UIOS property for (1.4), **without** knowledge of the exact solution map $\phi$, we must first select appropriate functions $g_i \in C^1(\Re^n; \Re)$ ($i = 1,\ldots,k$) and show that the following (more demanding) inequalities hold:

$$\nabla V_i(x) f(x, x_0, d, v, v_0) \leq -\rho_i(V_i(x)),$$

$$\forall (x,d,v,v_0) \in \Re^n \times D \times U \times U \text{ with } \max\{\zeta(|v|),\zeta(|v_0|)\} \leq V_i(x) \text{ and } x_0 \in B_i^{g_i}(h(x_0),x) \ (i=1,...,k) \quad (3.4)$$

instead of (3.2). The examples of the following section illustrate the use of Theorem 3.1 in conjunction with (3.4).

**Proof of Theorem 3.1:** Consider a solution $x(t)$ of (1.4) under hypotheses (H1-4) corresponding to arbitrary $(v,d,\tilde{d}) \in M(U) \times M(D) \times M(\Re^+)$ with initial condition $x(0) = x_0 \in \Re^n$. By virtue of Theorem 2.1, there exists a maximal existence time for the solution denoted by $t_{\max} \leq +\infty$. Let $V_i(t) = V_i(x(t))$, $i=1,...,k$, absolutely continuous functions on $[0,t_{\max})$. Moreover, let $V(t) := \max_{i=1,...,k} V_i(t)$, $\pi := \{\tau_0,\tau_1,...\}$ the set of sampling times (which may be finite if $t_{\max} < +\infty$) and $p(t) := \max\{\tau \in \pi : \tau \leq t\}$, $q(t) := \min\{\tau \in \pi : \tau \geq t\}$. Let $I \subset [0,t_{\max})$ be the zero Lebesgue measure set where $x(t)$ is not differentiable. Clearly, we have $x(t) = \phi(t - p(t), x(p(t)); P_t d, P_t v)$ for all $t \in [0,t_{\max})$, where $(P_t v)(s) = v(p(t)+s)$, $(P_t d)(s) = d(p(t)+s)$, $s \geq 0$. Next we show that the following implication holds for $t \in [0,t_{\max}) \setminus I$ and $i=1,...,k$:

$$V_i(t) \geq \max\left\{\sup_{p(t) \leq s \leq t} \zeta(|v(s)|), \sup_{p(t) \leq s \leq t} a(V(s))\right\} \Rightarrow \dot{V}_i(t) \leq -\rho_i(V_i(t)) \quad (3.5)$$

In order to prove implication (3.5) let $t \in [0,t_{\max}) \setminus I$, $i=1,...,k$, $\tau = p(t)$ and suppose that $V_i(t) \geq \max\left\{\sup_{p(t) \leq s \leq t} \zeta(|v(s)|), \sup_{p(t) \leq s \leq t} a(V(s))\right\}$. By virtue of the semigroup property for the previous inequality implies that $\zeta(|v(\tau+s)|) = \zeta(|(P_t v)(s)|) \leq V_i(x(t))$, $a(V_j(\phi(s,x(\tau);P_t d,P_t v))) \leq V_i(x(t))$ for all $s \in [0,t-\tau]$ and $j=1,...,k$. In this case, by virtue of definition (3.3) and the fact that $t-\tau \leq h(x(\tau))$, it follows that $x(\tau) \in A_i(h(x(\tau)),x(t))$. Since $\dot{x}(t) = f(x(t),x(\tau),d(t),v(t),v(\tau))$, we conclude from (3.2) that $\dot{V}_i(t) \leq -\rho_i(V_i(t))$.

Lemma 2.4 implies that there exists a family of continuous function $\sigma_i$ of class $KL$ ($i=1,...,k$), with $\sigma_i(s,0)=s$ for all $s \geq 0$ (which is independent of $(x_0,v,d,\tilde{d}) \in \Re^n \times M(U) \times M(D) \times M(\Re^+)$) such that for all $\xi \in [0,t_{\max})$, $t \in [\xi,t_{\max})$ and $i=1,...,k$ we have:

$$V_i(t) \leq \max\left\{\sigma_i(V_i(\xi),t-\xi); \sup_{\xi \leq \tau \leq t} \sigma_i\left(a\left(\sup_{p(\tau) \leq s \leq \tau} V(s)\right),t-\tau\right); \sup_{\xi \leq \tau \leq t} \sigma_i\left(\zeta\left(\sup_{p(\tau) \leq s \leq \tau} |v(s)|\right),t-\tau\right)\right\} \quad (3.6)$$

Let $\sigma(s,t) := \max_{i=1,...,k} \sigma_i(s,t)$, which is a function of class $KL$ that satisfies $\sigma(s,0)=s$ for all $s \geq 0$. An immediate consequence of the previous definition, estimate (3.6) with $\xi=0$ and the fact that $\sigma_i(s,0)=s$ for all $s \geq 0$ ($i=1,...,k$) is the following estimate, which holds for all $t \in [0,t_{\max})$:

$$V(t) \leq \max\left\{\sigma(V(0),t); a\left(\sup_{0 \leq s \leq t} V(s)\right); \sup_{0 \leq s \leq t} \zeta(|v(s)|)\right\} \quad (3.7)$$

Inequality (3.7) and the fact that $\sigma(s,0)=s$ for all $s \geq 0$, imply the following inequality:

$$\sup_{0 \leq s \leq t} V(s) \leq \max\left\{V(0); a\left(\sup_{0 \leq s \leq t} V(s)\right); \sup_{0 \leq s \leq t} \zeta(|v(s)|)\right\}, \text{ for all } t \in [0,t_{\max}) \quad (3.8)$$

Making use of the fact that $a(s) < s$ for all $s > 0$, we obtain from (3.8):

$$V(t) \leq \max\left\{V(0) \ ; \sup_{0 \leq s \leq t} \zeta(|v(s)|)\right\}, \text{ for all } t \in [0,t_{\max}) \quad (3.9)$$

Clearly, inequality (3.9) implies that as long as the solution of (1.4) exists, $V(t)$ is bounded. A standard contradiction argument in conjunction with the "Boundedness-Implies-Continuation" property for (1.4) (see Theorem 2.1), inequality (3.1) and hypothesis (H3) shows that necessarily we must have $t_{max} = +\infty$. We conclude that estimate (3.6) holds for all $\xi \geq 0$, $t \geq \xi$ and $i = 1,...,k$. Similarly, estimate (3.9) holds for all $t \geq 0$.

An immediate consequence of estimate (3.6) with $\xi \in \pi$ (where $\pi := \{\tau_0, \tau_1,...\}$ is the set of sampling times), estimate (3.9), definitions $V(t) := \max_{i=1,...,k} V_i(t)$, $\sigma(s,t) := \max_{i=1,...,k} \sigma_i(s,t)$ and the fact that $\sigma_i(s,0) = s$ for all $s \geq 0$ ($i = 1,...,k$) is the following estimate, which holds for all $t \geq \xi$, $\xi \in \pi$:

$$V(t) \leq \max\left\{\sigma(V(0), t-\xi); a\left(\sup_{\xi \leq s \leq t} V(s)\right); \sup_{0 \leq s \leq t} \zeta(|v(s)|)\right\} \quad (3.10)$$

We next show that the following estimate holds for all $\xi \geq 0$, $t \geq \xi$:

$$V(t) \leq \max\left\{\tilde{\sigma}(V(0), t-\xi); a\left(\sup_{\xi \leq s \leq t} V(s)\right); \sup_{0 \leq s \leq t} \zeta(|v(s)|)\right\} \quad (3.11)$$

where $\tilde{\sigma}(s,t) := s\exp(r-t)$ for $t \in [0,r]$ and $\tilde{\sigma}(s,t) := \sigma(s, t-r)$ for $t > r$, $r > 0$ being the upper bound for $h$. Since $\sigma(s,0) = s$ for all $s \geq 0$, it follows that $\tilde{\sigma}$ is of class $KL$.

First notice that estimate (3.10) implies estimate (3.11) when $\xi \in \pi$. We continue by showing that estimate (3.11) holds for the case $\xi \notin \pi$, $\xi \geq t_0$, as well. Let arbitrary $\xi \notin \pi$, $\xi \geq t_0$. Since $\xi \notin \pi$, it follows that $p(\xi) < \xi < q(\xi)$ (where $p(t) := \max\{\tau \in \pi : \tau \leq t\}$, $q(t) := \min\{\tau \in \pi : \tau \geq t\}$). Notice that by virtue of (3.9) and previous definition of $\tilde{\sigma} \in KL$, it follows that (3.11) holds for all $t \leq \xi + r$. Inequality (3.6) gives for all $t \geq \xi + r \geq q(\xi)$ and $i = 1,...,k$:

$$V_i(t) \leq \max\left\{\begin{array}{l}\sigma_i(V_i(\xi), t-\xi); \sup_{\xi \leq \tau < q(\xi)} \sigma_i\left(a\left(\sup_{p(\tau) \leq s \leq \tau} V(s)\right), t-\tau\right); \\ \sup_{q(\xi) \leq \tau \leq t} \sigma_i\left(a\left(\sup_{p(\tau) \leq s \leq \tau} V(s)\right), t-\tau\right); \sup_{p(\xi) \leq s \leq t} \zeta(|v(s)|)\end{array}\right\} \quad (3.12)$$

Notice that if $\xi \leq \tau < q(\xi)$ then $p(\tau) = p(\xi)$. Moreover, if $\tau \geq q(\xi)$ then $p(\tau) \geq \xi$. Thus, by making use of the facts that $a(s) \leq s$ for all $s \geq 0$ and $q(\xi) \leq \xi + r$, we obtain from (3.12):

$$V_i(t) \leq \max\left\{\begin{array}{l}\sigma_i(V_i(\xi), t-\xi); \sup_{\xi \leq \tau < \xi+r} \sigma_i\left(\sup_{p(\xi) \leq s \leq \tau} V(s), t-\tau\right); \\ a\left(\sup_{\xi \leq s \leq t} V(s)\right); \sup_{p(\xi) \leq s \leq t} \zeta(|v(s)|)\end{array}\right\}, \forall t \geq \xi + r, \ i = 1,...,k \quad (3.13)$$

Inequality (3.9) in conjunction with inequality (3.13) and the fact that $\sigma_i(s,0) = s$ for all $s \geq 0$ ($i = 1,...,k$), gives:

$$V_i(t) \leq \max\left\{\sigma_i(V(0), t-\xi-r); a\left(\sup_{\xi \leq s \leq t} V(s)\right); \sup_{0 \leq s \leq t} \zeta(|v(s)|)\right\}, \forall t \geq \xi + r, \ i = 1,...,k \quad (3.14)$$

Since $\sigma(s,t) := \max_{i=1,...,k} \sigma_i(s,t)$ and $V(t) = \max_{i=1,...,k} V_i(t)$, we obtain from (3.14):

$$V(t) \leq \max\left\{\sigma(V(0), t-\xi-r); a\left(\sup_{\xi \leq s \leq t} V(s)\right); \sup_{0 \leq s \leq t} \zeta(|v(s)|)\right\}, \forall t \geq \xi + r \quad (3.15)$$

We conclude from (3.15) and definition of $\tilde{\sigma} \in KL$, that (3.11) holds for that case $t \geq \xi + r$, as well.

Lemma 2.3 in conjunction with inequality (3.11) implies the existence of $\bar{\sigma} \in KL$ such that:

$$V(t) \leq \max\left\{ \bar{\sigma}(V(0),t);\ \sup_{0\leq s\leq t} \zeta(|v(s)|) \right\} \tag{3.16}$$

Clearly, estimate (3.16) in conjunction with inequality (3.1) shows that (1.4) satisfies the UIOS property with gain $\gamma = a_1^{-1} \circ \zeta \in N$ from the input $v \in M(U)$ and zero gain from the input $\tilde{d} \in M(\Re^+)$. The proof is complete. ◁

The following corollary is an immediate consequence of Theorem 3.1 and the fact that for every positive definite and radially unbounded function $V \in C^1(\Re^n; \Re^+)$ and for every positive definite function $W \in C^0(\Re^n; \Re^+)$, there exists a positive definite function $\rho \in C^0(\Re^+; \Re^+)$ such that $\rho(V(x)) \leq W(x)$, for all $x \in \Re^n$.

**Corollary 3.3 (Single Lyapunov Function Characterization of UISS):** *Consider system (1.4) under hypotheses (H1-4) and suppose that there exists a function $V \in C^1(\Re^n; \Re^+)$, functions $a_1, a_2 \in K_\infty$, $a, \zeta \in N$ with $a(s) < s$ for all $s > 0$ and a positive definite mapping $W \in C^0(\Re^n; \Re^+)$, such that the following inequalities hold:*

$$a_1(|x|) \leq V(x) \leq a_2(|x|), \quad \forall x \in \Re^n \tag{3.17}$$

$$\nabla V(x) f(x, x_0, d, v, v_0) \leq -W(x),$$
$$\forall (x, d, v, v_0) \in \Re^n \times D \times U \times U \text{ with } \max\{\zeta(|v|), \zeta(|v_0|)\} \leq V(x) \text{ and } x_0 \in A(h(x_0), x) \tag{3.18}$$

*where the set-valued map $\Re^+ \times \Re^n \ni (h,x) \to A(h,x) \subseteq \Re^n$ is defined by*

$$A(h,x) = \bigcup_{0\leq s\leq h} \left\{ \begin{array}{l} x_0 \in \Re^n : \exists (d,v) \in M(D) \times M(U) \text{ with } \phi(s, x_0; d, v) = x \\ \zeta(|v(t)|) \leq V(x) \text{ and } a(V(\phi(t, x_0; d, v))) \leq V(x), \text{ for all } t \in [0, s] \end{array} \right\} \tag{3.19}$$

*and $\phi(t, x_0; d, v)$ denotes the solution of $\dot{x}(t) = f(x(t), x_0, d(t), v(t), v(0))$ with initial condition $x(0) = x_0$ corresponding to $(d, v) \in M(D) \times M(U)$.*

*Then (1.4) satisfies the UISS property with gain $\gamma = a_1^{-1} \circ \zeta \in N$ from the input $v \in M(U)$ and zero gain from the input $\tilde{d} \in M(\Re^+)$.*

It should be noted that in practice (where the solution map is rarely known), Corollary 3.3 is used in conjunction with Remark 3.2.

## 4. Illustrating Examples

The following example illustrates the use of Theorem 3.1 and Corollary 3.3. Moreover, it shows that the use of vector Lyapunov functions provides flexibility, which is not easily obtained with the use of a single Lyapunov function.

**Example 4.1:** Consider the planar control system

$$\dot{x}_1 = -2x_1 - d_1 x_1^3 + x_2 \quad ; \quad \dot{x}_2 = d_2 x_2^2 - x_2^3 + u$$
$$x = (x_1, x_2)' \in \Re^2, d = (d_1, d_2) \in D := [\delta, \Delta] \times [-1,1], u \in \Re \tag{4.1}$$

We consider the closed-loop system (4.1) with the sampled feedback $u(t) = -2x_2(\tau_i) + v(t)$, $t \in [\tau_i, \tau_{i+1})$ and positive sampling rate, i.e., we consider the system:

$$\dot{x}_1(t) = -2x_1(t) - d_1(t)x_1^3(t) + x_2(t) \quad , \quad t \in [\tau_i, \tau_{i+1})$$
$$\dot{x}_2(t) = d_2(t)x_2^2(t) - x_2^3(t) - 2x_2(\tau_i) + v(t), \quad t \in [\tau_i, \tau_{i+1})$$
$$\tau_0 = t_0, \tau_{i+1} = \tau_i + r\exp(-\tilde{d}(\tau_i)) \tag{4.2}$$
$$x(t) = (x_1(t), x_2(t))' \in \Re^2$$
$$(d_1(t), d_2(t), \tilde{d}(t), v(t)) \in [\delta, \Delta] \times [-1,1] \times \Re^+ \times \Re$$

where the input $v(t) \in U = \Re$ quantifies the effect of control actuator errors. It should be noted that the feedback law $u(t) = -2x_2(\tau_i)$ is obtained from the discretization of the globally stabilizing continuous-time feedback $u(t) = -2x_2(t)$ (emulation procedure). This example was selected so that:

i) results concerning globally Lipschitz systems (see [13]) cannot be applied since system (4.2) is not globally Lipschitz,
ii) results concerning homogeneous systems (see [8]) cannot be applied since system (4.2) is not homogeneous,
iii) if the results contained in [14,15,48] were applied to system (4.2) then they would lead to **local** stability properties,
iv) if the results contained in [24,34,35,36] were applied to system (4.2) then they would lead to **semiglobal practical** stability properties.

We next show that there exists $r, K > 0$ and $\delta, \Delta \in \Re^+$ with $\delta \le \Delta$ such that system (4.2) satisfies the UISS property with gain $\gamma(s) = K s \in N$ from the input $v \in M(\Re)$ and zero gain from the input $\tilde{d} \in M(\Re^+)$. In order to show the advantages of the use of vector Lyapunov functions, we perform the analysis using a single and a vector Lyapunov function.

Notice that for every $\delta, \Delta \in \Re^+$ with $\delta \le \Delta$ and $r > 0$, system (4.2) with output $H(x) = x$ satisfies hypotheses (H1-4).

1st Method: Single Lyapunov function

Consider the function:

$$V(x) = \frac{1}{2}x_1^2 + \frac{1}{2}x_2^2 \tag{4.3}$$

Clearly, the function $V \in C^1(\Re^n; \Re^+)$ defined by (4.3) satisfies (3.17) with $a_1(s) = a_2(s) = \frac{1}{2}s^2$. Moreover, by virtue of (4.2) and definition (4.3), we obtain for all $(x_1, x_2, x_{2,0}, d_1, d_2, v) \in \Re^3 \times [\delta, \Delta] \times [-1,1] \times \Re$:

$$\nabla V(x) \begin{bmatrix} -2x_1 - d_1 x_1^3 + x_2 \\ d_2 x_2^2 - x_2^3 - 2x_{2,0} + v \end{bmatrix} = -2x_1^2 - \delta x_1^4 + x_1 x_2 + d_2 x_2^3 - x_2^4 - 2x_2 x_{2,0} + x_2 v$$
$$\le -\frac{3}{2}x_1^2 - \delta x_1^4 - x_2^2 - \frac{1}{2}x_2^4 - 2x_2(x_{2,0} - x_2) + |x_2||v| \tag{4.4}$$

Let $g(x) := x_2$, $\zeta(s) := 2s^2 \in N$ and $a(s) := c^{-2}s \in N$, where $c > 1$. It follows from (4.2) and (4.3) that the mapping $b^g(x) := \max\{|\nabla g(\xi)f(\xi, x_0, d, v)| : d \in D, a(\max\{V(\xi), V(x_0)\}) \le V(x), \zeta(|v|) \le V(x)\}$ satisfies the following inequality for all $(x_1, x_2) \in \Re^2$:

$$b^g(x) \le c^2(x_1^2 + x_2^2) + c^3(x_1^2 + x_2^2)^{\frac{3}{2}} + \left(2c + \frac{1}{2}\right)(x_1^2 + x_2^2)^{\frac{1}{2}}$$

Using Young inequalities and the above inequality, we obtain for all $(x_1, x_2) \in \Re^2$:

$$|x_2|b^g(x) \le \frac{5c^2+1}{2}x_2^2 + \left(c+\frac{1}{4}\right)x_1^2 + \frac{5c^3}{2}x_1^4 + \frac{11c^3}{4}x_2^4 \qquad (4.5)$$

Since $A(h,x) \subseteq B^g(h,x) = \{x_0 \in \Re^2 : |g(x_0) - g(x)| \le h b^g(x)\}$ for all $(h,x) \in \Re^+ \times \Re^2$ (see Remark 3.2), it follows from (4.4), (4.5) that the following inequality holds for all $(h,x,d) \in \Re^+ \times \Re^2 \times D$, $\zeta(|v|) \le V(x)$ and $(x_{1,0}, x_{2,0})' \in A(h,x)$:

$$\nabla V(x)\begin{bmatrix} -2x_1 - d_1 x_1^3 + x_2 \\ d_2 x_2^2 - x_2^3 - 2x_{2,0} + v \end{bmatrix} \le -\frac{1}{8}(7-16ch-4h)x_1^2 - (\delta-5c^3h)x_1^4 - \frac{1}{8}(7-40c^2h-8h)x_2^2 - \frac{1}{2}(1-11c^3h)x_2^4 \quad (4.6)$$

Inequality (4.6) guarantees that (3.18) holds if $\delta > 0$ and if we define $h(x) := r > 0$ with

$$r < \frac{7}{40c^2+8} \quad \text{and} \quad r \le c^{-3}\min\left\{\frac{1}{11}; \frac{\delta}{5}\right\} \qquad (4.7)$$

and $W(x) = \mu|x|^2$ for $\mu := \frac{1}{8}(7 - 40c^2 r - 8r)$. We conclude from Corollary 3.3 that that system (4.2) satisfies the UISS property with gain $\gamma(s) = Ks \in N$ (with $K = 2$) from the input $v \in M(\Re)$ and zero gain from the input $\tilde{d} \in M(\Re^+)$.

2nd Method: Vector Lyapunov function

Consider the functions:

$$V_1(x) = \frac{1}{2}x_1^2; \quad V_2(x) = \frac{1}{2}x_2^2 \qquad (4.8)$$

Clearly, the functions $V_i \in C^1(\Re^n; \Re^+)$, $i=1,2$ defined by (4.8) satisfy (3.1) with $H(x) := x$, $a_1(s) = \frac{1}{4}s^2$, $a_2(s) = \frac{1}{2}s^2$ (an immediate consequence of the inequality $\frac{1}{2}(a+b) \le \max\{a,b\} \le a+b$, which holds for all $a, b \ge 0$). Let $\zeta(s) := 2s^2 \in N$ and $a(s) := c^{-2}s \in N$, where $c \in (1,2)$.

By virtue of (4.2) and definition (4.8), we obtain for all $\delta \ge 0$, $(x_1, x_2, d_1) \in \Re^2 \times [\delta, \Delta]$ with $a(V_2(x)) \le V_1(x)$:

$$\nabla V_1(x)(-2x_1 - d_1 x_1^3 + x_2) = -2x_1^2 - \delta x_1^4 + x_1 x_2 \le -(2-c)x_1^2 = -\frac{2-c}{2}V_1(x)$$

The above inequality guarantees that (3.2) holds for $i=1$ with $\rho_1(s) = \frac{2-c}{2}s$.

Furthermore, we obtain from (4.2) and definition (4.8), for all $(x_1, x_2, d_2, v) \in \Re^2 \times [-1,1] \times \Re$ with $\zeta(|v|) \le V_2(x)$:

$$\nabla V_2(x)(d_2 x_2^2 - x_2^3 - 2x_{2,0} + v) = d_2 x_2^3 - x_2^4 - 2x_2 x_{2,0} + x_2 v \le -x_2^2 - \frac{1}{2}x_2^4 - 2x_2(x_{2,0} - x_2) \qquad (4.9)$$

Let $g(x) := x_2$. It follows from (4.2) and (4.8) that the mapping $b^g(x) := \max\{|\nabla g(\xi)f(\xi, x_0, d, v)| : d \in D, a(\max\{V(\xi), V(x_0)\}) \le V_2(x), \zeta(|v|) \le V_2(x)\}$ satisfies the following inequality for all $(x_1, x_2) \in \Re^2$:

$$b_2^g(x) \le c^2 x_2^2 + c^3 |x_2|^3 + \left(2c + \frac{1}{2}\right)|x_2|$$

Since $A_2(h,x) \subseteq B_2^g(h,x) = \{x_0 \in \Re^2 : |g(x_0) - g(x)| \le h\, b_2^g(x)\}$ for all $(h,x) \in \Re^+ \times \Re^2$ (see Remark 3.2), it follows from the above inequality and (4.9) that the following inequality holds for all $(x_1, x_2, d_2, v) \in \Re^2 \times [-1,1] \times \Re$ with $\zeta(|v|) \le V_2(x)$ and $(x_{1,0}, x_{2,0})' \in A_2(h,x)$:

$$\nabla V_2(x)\left(d_2 x_2^2 - x_2^3 - 2x_{2,0} + v\right) \le -(1 - 2h - 5c^2 h)x_2^2 - \frac{1}{2}(1 - 6c^3 h)x_2^4 \qquad (4.10)$$

Inequality (4.10) guarantees that (3.2) for $i = 2$ holds if we define $h(x) := r > 0$ with

$$r < \frac{1}{5c^2 + 2} \text{ and } r \le \frac{1}{6c^3} \qquad (4.11)$$

and $\rho_2(s) = \dfrac{1 - 2r - 5c^2 r}{2} s$. We conclude from Theorem 3.1 that that system (4.2) satisfies the UISS property with gain $\gamma(s) = K s \in \mathcal{N}$ (with $K = 2\sqrt{2}$) from the input $v \in M(\Re)$ and zero gain from the input $\tilde{d} \in M(\Re^+)$.

Comparison of the two methods

We are now ready to compare the two methods. It is clear from (4.7) and (4.11) that the vector Lyapunov function method is less conservative compared to the single Lyapunov method since:

i) it allows greater sampling periods
ii) we do not have to assume that $\delta > 0$ (only $\delta \ge 0$ suffices and notice that the size of $\delta$ does not affect the sampling rate)

Of course, different Lyapunov functions may be used in order to obtain even less conservative results. ◁

Example 4.1 shows that there exist nonlinear systems with no special characteristics (such as homogeneity or global Lipschitz properties), which can be globally asymptotically stabilized by means of sampled-data feedback with positive sampling rate. Working exactly as in the previous example with a vector Lyapunov function, we can establish that $0 \in \Re^2$ can be robustly globally asymptotically stabilized by linear feedback with zero order hold and positive sampling rate for the system:

$$\begin{aligned}\dot{x}_1 &= f_1(x_1, x_2, d) \\ \dot{x}_2 &= f_2(x_1, x_2, d) + u \\ x &= (x_1, x_2)' \in \Re^2, d \in D, u \in \Re\end{aligned} \qquad (4.12)$$

where $D \subset \Re^m$ is a compact set, $f_i : \Re^2 \times D \to \Re$ ($i = 1,2$) are locally Lipschitz mappings with $f_i(0,0,d) = 0$ for all $d \in D$ ($i = 1,2$), under the following hypothesis:

**(P)** *There exist constants* $c > 1$, $a \in \Re$, $L, \gamma \ge 0$ *such that:*

$$\max\{x_1 f_1(x_1, -ax_1 + \xi, d) ; x_1 \ne 0, |\xi| \le c|x_1|, d \in D\} < 0, \text{ for all } x_1 \ne 0 \qquad (4.13a)$$

$$\begin{aligned}&|z|\max\{|f_2(x_1, -ax_1 + \xi, d) + a f_1(x_1, -ax_1 + \xi, d)|; \max(|x_1|, |\xi|) \le c|z|, d \in D\} \\ &+ L \max\{z f_2(x_1, -ax_1 + z, d) + a z f_1(x_1, -ax_1 + z, d) ; |x_1| \le c|z|, d \in D\} \le \gamma z^2\end{aligned}, \text{ for all } z \in \Re \qquad (4.13b)$$

Indeed, under hypothesis (P) we may conclude that there exist constants $r, R > 0$ such that the following system:

$$\dot{x}_1(t) = f_1(x_1(t), x_2(t), d(t)), \, t \in [\tau_i, \tau_{i+1})$$
$$\dot{x}_2(t) = f_2(x_1(t), x_2(t), d(t)) - R(x_2(\tau_i) + ax_1(\tau_i)), \, t \in [\tau_i, \tau_{i+1}) \quad (4.14)$$
$$\tau_{i+1} = \tau_i + r \exp(-\tilde{d}(\tau_i))$$

where $a \in \Re$ is the constant involved in (4.13), satisfies the UISS property with zero gain from the input $\tilde{d} \in M(\Re^+)$. For example, it can be verified that system (4.1) satisfies hypothesis (P). The vector Lyapunov function that can be used for the stability analysis is defined by the equations $V_1(x) = \frac{1}{2}x_1^2$ and $V_2(x) = \frac{1}{2}(x_2 + ax_1)^2$, where $a \in \Re$ is the constant involved in hypothesis (P).

**Example 4.2:** We consider triangular nonlinear systems of the form:

$$\dot{x}_i = \sum_{j=1}^{i} x_j \phi_{i,j}(x_1, \ldots, x_i, d) + g_i(x_1, \ldots, x_i, d) x_{i+1} \quad , \quad i = 1, \ldots, n-1$$
$$\dot{x}_n = \sum_{j=1}^{n} x_j \phi_{n,j}(x_1, \ldots, x_n, d) + g_n(x_1, \ldots, x_n, d) u \quad (4.15)$$
$$x = (x_1, \ldots, x_n)' \in \Re^n, d \in D, u \in \Re$$

where $D \subset \Re^m$ is a compact set, $\phi_{i,j} : \Re^i \times D \to \Re$ ($j = 1, \ldots, i$, $i = 1, \ldots, n$), $g_i : \Re^i \times D \to \Re$ ($i = 1, \ldots, n$) locally Lipschitz mappings with $g_i(x_1, \ldots, x_i, d) > 0$ for all $(x, d) \in \Re^n \times D$ ($i = 1, \ldots, n$). A recursive method (backstepping) is presented in [4] (but see also the references therein) for the construction of a positive definite and radially unbounded function $V \in C^1(\Re^n; \Re^+)$, a mapping $\zeta \in N$, a positive definite function $W \in C^1(\Re^n; \Re^+)$ and a feedback function $k \in C^0(\Re^n; \Re)$ with $k(0) = 0$, such that:

$$\zeta(|e|) \leq V(x) \Rightarrow \nabla V(x)(F(x,d) + G(x,d)k(x+e)) \leq -W(x), \, \forall x \in \Re^n \quad (4.16a)$$
(guarantees robustness with respect to additive measurement errors)

or

$$\zeta(|v|) \leq V(x) \Rightarrow \nabla V(x)(F(x,d) + G(x,d)k(x) + G(x,d)v) \leq -W(x), \, \forall x \in \Re^n \quad (4.16b)$$
(guarantees robustness with respect to control actuator errors)

where

$$F(x,d) := \left( x_1 \phi_{1,1}(x_1,d) + g_1(x_1,d)x_2, \ldots, \sum_{j=1}^{n-1} x_j \phi_{n-1,j}(x_1, \ldots, x_{n-1}, d) + g_{n-1}(x_1, \ldots, x_{n-1}, d)x_n, \sum_{j=1}^{n} x_j \phi_{n,j}(x_1, \ldots, x_n, d) \right)',$$

$G(x,d) := (0, \ldots, 0, g_n(x,d))'$. We consider the question of whether the constructed feedback $k \in C^0(\Re^n; \Re)$ can robustly globally asymptotically stabilize $0 \in \Re^n$ for system (4.15) when applied with zero order hold and positive sampling rate as well as the determination of the maximum allowable sampling period for which robust global asymptotic stability is preserved.

Corollary 3.3 provides a direct answer to both issues stated above. Particularly, we have:

If there exists $a \in N$ with $a(s) < s$ for all $s > 0$ and a constant $h > 0$ such that:

$$\zeta(h\rho(x)) \leq V(x), \text{ for all } x \in \Re^n \quad (4.17)$$

where

$$\rho(x) := \max\{|F(\xi,d) + G(\xi,d)k(x_0)| : d \in D, a(\max\{V(\xi); V(x_0)\}) \leq V(x)\} \quad (4.18a)$$

in case that (4.16a) holds or

$$\rho(x) := \max\{|\nabla k(\xi)(F(\xi,d)+G(\xi,d)k(x_0))|: d \in D, a(\max\{V(\xi); V(x_0)\}) \le V(x)\} \quad (4.18b)$$

in case that (4.16b) holds and $k \in C^1(\Re^n; \Re)$, then the constructed feedback $k \in C^0(\Re^n; \Re)$ can robustly globally asymptotically stabilize $0 \in \Re^n$ for system (4.15) when applied with zero order hold and positive sampling rate. Moreover, the maximum allowable sampling period is the constant $h > 0$ involved in (4.17). Specifically, the sampled-data system:

$$\dot{x}_i(t) = \sum_{j=1}^{i} x_j \phi_{i,j}(x_1,...,x_i,d) + g_i(x_1(t),...,x_i(t),d(t))x_{i+1}(t) \quad , \quad i=1,...,n-1, \ t \in [\tau_i, \tau_{i+1})$$

$$\dot{x}_n(t) = \sum_{j=1}^{n} x_j \phi_{n,j}(x_1,...,x_i,d) + g_n(x_1(t),...,x_n(t),d(t))k(x(\tau_i)) \quad , \quad t \in [\tau_i, \tau_{i+1}) \quad (4.19)$$

$$\tau_{i+1} = \tau_i + h\exp(-\tilde{d}(\tau_i))$$

satisfies the UISS property with zero gain from the input $\tilde{d} \in M(\Re^+)$. The conclusion is important, since conditions (4.17), (4.18) can be easily checked after the application of the backstepping method presented in [4]. It should be emphasized that experimentation with the backstepping technique showed that it is more likely to guarantee stability for (4.19) when the feedback law $k \in C^0(\Re^n; \Re)$ is designed to ensure robustness with respect to additive state measurement errors (i.e., following the procedure described in Chapter 6 in [4]).   ◁

## 5. Concluding Remarks

In the present work sufficient conditions expressed by means of single and vector Lyapunov functions of Uniform Input-to-Output Stability (UIOS) and Uniform Input-to-State Stability (UISS) are given for finite-dimensional systems under feedback control with zero order hold (Corollary 3.3 and Theorem 3.1, respectively). Illustrating examples are presented, which show the flexibility that a vector Lyapunov function can provide. Moreover, we show how the main results can be used in conjunction with the backstepping method for triangular nonlinear systems to check whether the feedback, constructed by the backstepping method, can robustly globally asymptotically stabilize the equilibrium point when applied with zero order hold and positive sampling rate (emulation procedure).

It is clear that the main results of this work may be extended to the time-varying case by using time-varying Lyapunov or Razumikhin functions for time-delay systems. Moreover, the results can be extended to cover non-uniform and weighted notions of IOS and ISS (see [18,19]). Both extensions are straightforward in the sense that the proofs follow the same procedure as the proof of Theorem 3.1.

It should be emphasized that the sufficient conditions for ISS and IOS for sampled-data systems presented in this work are not restricted to sampled-data feedback laws derived by emulation. Theorem 3.1 and Corollary 3.3 can be used for the Lyapunov redesign procedure, as explained in [37]. This topic is very important, since very rarely sampled-data feedback designs produced by emulation can guarantee global asymptotic stability.

## References

[1] Artstein, Z. and G. Weiss, "State Nullification by Memoryless Output Feedback", *Mathematics of Control, Signals and Systems*, 17, 2005, 38-56.
[2] Bernstein, D.S. and C.V. Hollot, "Robust Stability for Sampled-Data Control Systems", *Systems and Control Letters*, 13, 1989, 217-226.
[3] Clarke, F.H., Y.S. Ledyaev, E.D. Sontag and A.I. Subbotin, "Asymptotic Controllability implies Feedback Stabilization", *IEEE Transactions on Automatic Control*, 42(10), 1997, 1394-1407.
[4] Freeman, R. A. and P. V. Kokotovic, "Robust Nonlinear Control Design- State Space and Lyapunov Techniques", Birkhauser, Boston, 1996.
[5] Fridman, E., A. Seuret and J.-P. Richard, "Robust Sampled-Data Stabilization of Linear Systems: An Input-Delay Approach", *Automatica*, 40, 2004, 1441-1446.

[6] Fridman, E., U. Shaked and V. Suplin, "Input/Output Delay Approach to Robust Sampled-Data $H_\infty$ Control", *Systems and Control Letters*, 54, 2005, 271-282.

[7] Grune, L., "Stabilization by Sampled and Discrete Feedback with Positive Sampling Rate", in "Stability and Stabilization of Nonlinear Systems", D. Aeyels, F. Lamnabhi-Lagarrigue and A. van der Schaft (Eds), Springer-Verlag, London, 1999, 165-182.

[8] Grune, L., "Homogeneous State Feedback Stabilization of Homogenous Systems", *SIAM Journal on Control and Optimization*, 38(4), 2000, 1288-1308.

[9] Grune, L., "Input-to-State Dynamical Stability and its Lyapunov Function Characterization", *IEEE Transactions on Automatic Control*, 47, 2002, 1499-1504.

[10] Grune, L. and D. Nesic, "Optimization Based Stabilization of Sampled-Data Nonlinear Systems via Their Approximate Discrete-Time Models", *SIAM Journal on Control and Optimization*, 42, 2003, 98-122.

[11] Grune, L., D. Nesic and J. Pannek, "Model Predictive Control for Nonlinear Sampled-Data Systems", in Assessment and Future Directions of Nonlinear Model Predictive Control (NMPC 05), F. Allgower, L. Biegler and R. Findeisen, Lecture Notes in Control and Information Sciences, in print.

[12] Hale, J.K. and S.M.V. Lunel, "Introduction to Functional Differential Equations", Springer-Verlag, New York, 1993.

[13] Herrmann, G., S.K. Spurgeon and C. Edwards, "Discretization of Sliding Mode Based Control Schemes", *Proceedings of the 38th Conference on Decision and Control*, Phoenix, Arizona, U.S.A., 1999, 4257-4262.

[14] Hu, B. and A.N. Michel, "Stability Analysis of Digital Control Systems with time-varying sampling periods", *Automatica*, 36, 2000, 897-905.

[15] Hu, B. and A.N. Michel, "Robustness Analysis of Digital Control Systems with time-varying sampling periods", *Journal of the Franklin Institute*, 337, 2000, 117-130.

[16] Jiang, Z.P., A. Teel and L. Praly, "Small-Gain Theorem for ISS Systems and Applications", *Mathematics of Control, Signals and Systems*, 7, 1994, 95-120.

[17] Karafyllis, I., "The Non-Uniform in Time Small-Gain Theorem for a Wide Class of Control Systems with Outputs", *European Journal of Control*, 10(4), 2004, 307-323.

[18] Karafyllis, I., "A System-Theoretic Framework for a Wide Class of Systems I: Applications to Numerical Analysis", to appear in the *Journal of Mathematical Analysis and Applications*. Available in electronic form in the web site of the *Journal of Mathematical Analysis and Applications*.

[19] Karafyllis, I., "A System-Theoretic Framework for a Wide Class of Systems II: Input-to-Output Stability", to appear in the *Journal of Mathematical Analysis and Applications*. Available in electronic form in the web site of the *Journal of Mathematical Analysis and Applications*.

[20] Kazantzis, N. and C. Kravaris, "System-Theoretic Properties of Sampled-Data Representations of Nonlinear Systems Obtained via Taylor-Lie Series", *International Journal of Control*, 67(6), 1997, 997-1020.

[21] Kellett, C.M., H. Shim and A.R. Teel, "Further Results on Robustness of (possibly discontinuous) Sample and Hold Feedback", *IEEE Transactions on Automatic Control*, 49(7), 2004, 1081-1089.

[22] Khalil, H.K., "Nonlinear Systems", 2nd Edition, Prentice-Hall, 1996.

[23] Khalil, H.K., "Performance Recovery Under Output Feedback Sampled-Data Stabilization of a Class of Nonlinear Systems", *IEEE Transactions on Automatic Control*, 49(12), 2004, 2173-2184.

[24] Laila, D.S., D. Nesic and A.R. Teel, "Open and Closed-Loop Dissipation Inequalities Under Sampling and Controller Emulation", *European Journal of Control*, 18, 2002, 109-125.

[25] Laila, D. S. and A. Astolfi, "Input-to-State Stability for Discrete-time Time-Varying Systems with Applications to Robust Stbilization of Systems in Power Form", *Automatica*, 41, 2005, 1891-1903.

[26] Lin, Y., E.D. Sontag and Y. Wang, "A Smooth Converse Lyapunov Theorem for Robust Stability", *SIAM Journal on Control and Optimization*, 34, 1996, 124-160.

[27] Mancilla-Aguilar, J.L., R. A. Garcia and M.I. Troparevsky, "Stability of a Certain Class of Hybrid Dynamical Systems", *International Journal of Control*, 73(15), 2000, 1362-1374.

[28] Monaco, S. and D. Normand-Cyrot, "Issues on Nonlinear Digital Control", *European Journal of Control*, 7, 2001, 160-177.

[29] Naghshtabrizi, P., J. P. Hespanha and A.R. Teel, "On the Robust Stability and Stabilization of Sampled-Data Systems: A Hybrid System Approach", to appear in the *Proceedings of 45th IEEE Conference on Decision and Control*, available electronically in http://www.ccec.ece.ucsb.edu/~payam/publication.htm.

[30] Nersesov, S.G. and W. M. Haddad, "On the Stability and Control of Nonlinear Dynamical Systems via Vector Lyapunov Functions", *IEEE Transactions on Automatic Control*, 51(2), 2006, 203-215.

[31] Nesic, D., A.R. Teel and E.D. Sontag, "Formulas Relating KL Stability Estimates of Discrete-Time and Sampled-Data Nonlinear Systems", *Systems and Control Letters*, 38(1), 1999, 49-60.

[32] Nesic, D., A.R. Teel and P.V. Kokotovic, "Sufficient Conditions for Stabilization of Sampled-Data Nonlinear Systems via Discrete-Time Approximations", *Systems and Control Letters*, 38(4-5), 1999, 259-270.

[33] Nesic, D. and A.R. Teel, "Sampled-Data Control of Nonlinear Systems: An Overview of Recent Results", Perspectives on Robust Control, R.S.O. Moheimani (Ed.), Springer-Verlag: New York, 2001, 221-239.

[34] Nesic, D. and D.S. Laila, "A Note on Input-to-State Stabilization for Nonlinear Sampled-Data Systems", *IEEE Transactions on Automatic Control*, 47, 2002, 1153-1158.

[35] Nesic, D. and D. Angeli, "Integral Versions of ISS for Sampled-Data Nonlinear Systems via Their Approximate Discrete-Time Models", *IEEE Transactions on Automatic Control*, 47, 2002, 2033-2038.

[36] Nesic, D. and A. Teel, "A Framework for Stabilization of Nonlinear Sampled-Data Systems Based on their Approximate Discrete-Time Models", *IEEE Transactions on Automatic Control*, 49(7), 2004, 1103-1122.

[37] Nesic, D. and L. Grune, "Lyapunov-based Continuous-time Nonlinear Controller Redesign for Sampled-Data Implementation", *Automatica*, 41, 2005, 1143-1156.

[38] Niculescu, S.I., "Delay Effects on Stability, A Robust Control Approach", Heidelberg, Germany, Springer-Verlag, 2001.

[39] Sontag, E.D., "Smooth Stabilization Implies Coprime Factorization", *IEEE Transactions on Automatic Control*, 34, 1989, 435-443.

[40] Sontag, E.D. and Y. Wang, "On Characterizations of the Input-to-State Stability Property", *Systems and Control Letters*, 24, 1995, 351-359.

[41] Sontag, E.D. and Y. Wang, "New Characterizations of the Input-to-State Stability", *IEEE Transactions on Automatic Control*, 41, 1996, 1283-1294.

[42] Sontag,, E.D., "Mathematical Control Theory", 2$^{nd}$ Edition, Springer-Verlag, New York, 1998.

[43] Sontag, E.D., "Clocks and Insensitivity to Small Measurement Errors", *ESAIM: Control, Optimisation and Calculus of Variations*, 4, 1999, 537-557.

[44] Sontag, E.D. and Y. Wang, "Notions of Input to Output Stability", *Systems and Control Letters*, 38, 1999, 235-248.

[45] Sontag, E.D. and Y. Wang, "Lyapunov Characterizations of Input-to-Output Stability", *SIAM Journal on Control and Optimization*, 39, 2001, 226-249.

[46] Sontag, E.D. and B. Ingalls, "A Small-Gain Theorem with Applications to Input/Output Systems, Incremental Stability, Detectability, and Interconnections", *Journal of the Franklin Institute*, 339, 2002, 211-229.

[47] Teel, A.R., "Connections between Razumikhin-Type Theorems and the ISS Nonlinear Small Gain Theorem", *IEEE Transactions on Automatic Control*, 43(7), 1998, 960-964.

[48] Ye, H., A.N. Michel and L. Hou, "Stability Theory for Hybrid Dynamical Systems", *IEEE Transactions on Automatic Control*, 43(4), 1998, 461-474.

[49] Zaccarian, L., A.R. Teel and D. Nesic, "On Finite Gain $L_p$ Stability of Nonlinear Sample-Data Systems", *Systems and Control Letters*, 49, 2003, 201-212.

# Appendix

**Proof of Lemma 2.4:** Notice that by virtue of Lemma 4.4 in [26], for each positive definite continuous function $\rho : \Re^+ \to \Re^+$ there exists a continuous function $\sigma$ of class $KL$, with $\sigma(s,0) = s$ for all $s \geq 0$ with the following property: if $y : [t_0, t_1] \to \Re^+$ is an absolutely continuous function and $I \subset [t_0, t_1]$ a set of Lebesgue measure zero such that $\dot{y}(t)$ is defined on $[t_0, t_1] \setminus I$ and such that the following differential inequality holds for all $t \in [t_0, t_1] \setminus I$:

$$\dot{y}(t) \leq -\rho(y(t)) \tag{A1}$$

then the following estimate holds for all $t \in [t_0, t_1]$:

$$y(t) \leq \sigma(y(t_0), t - t_0) \tag{A2}$$

Actually, the statement of Lemma 4.4 in [26] does not guarantee that $\sigma$ is continuous or that $\sigma(s,0) = s$ for all $s \geq 0$, but a close look at the proof of Lemma 4.4 in [26] shows that this is the case when $\rho : \Re^+ \to \Re^+$ is a positive definite continuous function. Moreover, notice that we may continuously extend $\sigma$ by defining $\sigma(s,t) := s \exp(-t)$ for $t < 0$.

Clearly, (2.8) holds for $t = t_0$ and $\sigma$ the function involved in (A2). We next show that (2.8) holds for arbitrary $t \in (t_0, t_1]$.

Let arbitrary $t \in (t_0, t_1]$ and define the functions $\tilde{u}(\tau) = \begin{cases} u(\tau) \text{ for } \tau \in [t_0, t] \\ 0 \text{ if otherwise} \end{cases}$, $\bar{u}(\tau) := \limsup_{\xi \to \tau} \tilde{u}(\xi)$. Notice that $\bar{u}$ is upper semi-continuous on $\tau \in [t_0, t]$ and consequently the function $p(\tau) := y(\tau) - \bar{u}(\tau)$ is lower semi-continuous on $[t_0, t]$. Next define the set:

$$A := \{\tau \in [t_0, t]: y(\tau) \leq \bar{u}(\tau)\} \tag{A3}$$

We distinguish the following cases:

1) $A = \emptyset$. In this case we have $y(\tau) > \bar{u}(\tau)$ for all $\tau \in [t_0, t]$. Since $\bar{u}(\tau) \geq u(\tau)$ for all $\tau \in [t_0, t]$, the previous inequality in conjunction with (2.7) implies that $\dot{y}(\tau) \leq -\rho(y(\tau))$ for all $\tau \in [t_0, t] \setminus I$. Thus in this case Lemma 4.4 in [26] guarantees that estimate (A2) holds.

2) $A \neq \emptyset$ and $\xi := \sup A < t$. In this case there exists a sequence $\tau_i \leq \xi$ with $\tau_i \to \xi$ and $y(\tau_i) - \bar{u}(\tau_i) \leq 0$. Since the function $p(t) = y(t) - \bar{u}(t)$ is lower semi-continuous, we obtain $p(\xi) = \liminf_{\tau \to \xi} p(\tau) \leq 0$ and consequently $y(\xi) \leq \bar{u}(\xi)$. Moreover, notice that by virtue of definition (A3), (2.7) and since $\bar{u}(\tau) \geq u(\tau)$ for all $\tau \in [t_0, t]$, the differential inequality $\dot{y}(\tau) \leq -\rho(y(\tau))$ holds for all $\tau \in (\xi, t] \setminus N$. Consequently, Lemma 4.4 in [26] implies $y(t) \leq \sigma(y(\tau), t - \tau)$ for all $\tau \in (\xi, t]$. By virtue of continuity of $\sigma$ and $y$ we get

$$y(t) \leq \sigma(y(\xi), t - \xi)$$

which combined $y(\xi) \leq \bar{u}(\xi)$ directly implies

$$y(t) \leq \sigma(\bar{u}(\xi), t - \xi) \leq \sup_{t_0 \leq s \leq t} \sigma(\bar{u}(s), t - s) \tag{A4}$$

3) $A \neq \emptyset$ and $\xi := \sup A = t$. In this case there exists a sequence $\tau_i \leq t$ with $\tau_i \to t$ and $y(\tau_i) - \bar{u}(\tau_i) \leq 0$. Since the function $p(t) = y(t) - \bar{u}(t)$ is lower semi-continuous, we obtain $p(t) = \liminf_{\tau \to t} p(\tau) \leq 0$ and consequently $y(t) \leq \bar{u}(t)$. Moreover, since $\sigma(s, 0) = s$ for all $s \geq 0$, it holds that $y(t) \leq u(t) = \sigma(\bar{u}(t), 0) \leq \sup_{t_0 \leq s \leq t} \sigma(\bar{u}(s), t - s)$.

Combining all the above cases, we may conclude that

$$y(t) \leq \max\left\{\sigma(y(t_0), t - t_0), \sup_{t_0 \leq s \leq t} \sigma(\bar{u}(s), t - s)\right\} \tag{A5}$$

Let $M := \sup_{t_0 \leq s \leq t} u(s)$. For each $\varepsilon > 0$, there exists $\delta > 0$ such that $\sigma(s, \tau - \delta) - \sigma(s, \tau) < \varepsilon$ for all $(s, \tau) \in [0, M] \times [0, t]$. Notice that since $\bar{u}(\tau) := \limsup_{\xi \to \tau} \tilde{u}(\xi)$ and $\tilde{u}(\tau) = \begin{cases} u(\tau) \text{ for } \tau \in [t_0, t] \\ 0 \text{ if otherwise} \end{cases}$, it follows that $\bar{u}(s) \leq \sup\{u(r): \max(s - \delta, t_0) \leq r \leq \min(s + \delta, t)\}$, for all $s \in [t_0, t]$. The previous inequalities imply:

$$\begin{aligned}\sigma(\bar{u}(s), t - s) &\leq \sup\{\sigma(u(r), t - s): \max(s - \delta, t_0) \leq r \leq \min(s + \delta, t)\} \\ &\leq \sup\{\sigma(u(r), t - r - \delta): \max(s - \delta, t_0) \leq r \leq \min(s + \delta, t)\} \\ &\leq \sup\{\sigma(u(r), t - r): \max(s - \delta, t_0) \leq r \leq \min(s + \delta, t)\} + \varepsilon \\ &\leq \sup_{t_0 \leq r \leq t} \sigma(u(r), t - r) + \varepsilon\end{aligned}$$

The above inequality in conjunction with (A5) imply that for each $\varepsilon > 0$, it holds that:

$$y(t) \leq \max\left\{\sigma(y(t_0), t - t_0), \sup_{t_0 \leq s \leq t} \sigma(u(s), t - s)\right\} + \varepsilon$$

Since $\varepsilon > 0$ is arbitrary, we conclude that the above estimate directly implies (2.8). The proof is complete. ◁